\documentclass[12pt,reqno]{article}

\usepackage[usenames]{color}
\usepackage{amssymb}
\usepackage{graphicx}
\usepackage{amscd}

\usepackage{amsthm}
\newtheorem{theorem}{Theorem}

\newtheorem{proposition}[theorem]{Proposition}

\theoremstyle{definition}

\newtheorem{example}[theorem]{Example}

\usepackage[colorlinks=true,
linkcolor=webgreen, filecolor=webbrown,
citecolor=webgreen]{hyperref}

\definecolor{webgreen}{rgb}{0,.5,0}
\definecolor{webbrown}{rgb}{.6,0,0}

\usepackage{color}

\usepackage{float}

\usepackage{graphics,amsmath,amssymb}
\usepackage{amsfonts}
\usepackage{latexsym}
\usepackage{epsf}

\begin{document}

\begin{center}
\vskip 1cm{\LARGE\bf Embedding structures associated with Riordan arrays and moment matrices} \vskip 1cm \large
Paul Barry\\
School of Science\\
Waterford Institute of Technology\\
Ireland\\
\href{mailto:pbarry@wit.ie}{\tt pbarry@wit.ie} \\

\end{center}
\vskip .2 in

\begin{abstract} Every ordinary Riordan array contains two naturally embedded Riordan arrays. We explore this phenomenon, and we compare it to the situation for certain moment matrices of families of orthogonal polynomials.
\end{abstract}

\section{Introduction} Riordan arrays \cite{SGWW} have been used mainly to prove combinatorial identities \cite{CMS, Spru}. Recently, their links to orthogonal polynomials have been investigated \cite{Barry_Meixner, Barry_Moment},  while there is a growing literature surrounding their structural properties \cite{Cheon, He, Jin, Alter}. In this note we investigate an embedding structure, common to all ordinary Riordan arrays. We also look at this embedding structure in the context of moment matrices of families of orthogonal polynomials. In addition to some knowledge of Riordan arrays, we assume that the reader has a basic familiarity with the theory of orthogonal polynomials on the real line \cite{Chihara, Gautschi, Szego}, production matrices \cite{Prod1, Prod2}, and continued fractions \cite{Wall}. We shall meet a number of integer sequences and integer triangles in this note. The On-Line Encyclopedia may be consulted for many of them \cite{SL1, SL2}.
In this note we shall understand by an ordinary Riordan array an integer number triangle whose $(n,k)$-th element $T_{n,k}$ is defined by a pair of power series $g(x)$ and $f(x)$ over the integers with $g(x)=1+g_1 x+g_2 x^2+\cdots$, $f(x)=x+f_2 x^2+f_3 x^3+\cdots$, in the following manner:
$$T_{n,k}=[x^n] g(x)f(x)^k.$$ 
The group law for Riordan arrays is given by \begin{displaymath}
(g,
f)\cdot (h,
l)=(g(h\circ f), l\circ f).\end{displaymath} 
If a matrix $A$ is the inverse of the coefficient array of a family of orthogonal polynomials, then we shall call it a moment matrix, and we shall single out the first column as the moment sequence. 

\section{The canonical embedding}
Let $(g, f)$ be an ordinary Riordan array $R$, with general term
$$T_{n,k}=[x^n] g f^k.$$

Then we observe that there are two naturally associated Riordan arrays ``embedded" in the array $R$ as follows.

Beginning at the first column of $R$, we take every second column, ``raising" the columns appropriately to obtain a lower-triangular matrix $A$.
The matrix $A$ is then the Riordan array
$$A=\left(g, \frac{f^2}{x}\right)$$ with general term $A_{n,k}$ given by
\begin{eqnarray*}
A_{n,k}&=& [x^n] g\left(\frac{f^2}{x}\right)^k \\
&=& [x^n] g x^{-k} f^{2k} \\
&=& [x^{n+k}] g f^{2k} \\
&=& T_{n+k,2k}.\end{eqnarray*}

\noindent Similarly, starting at the second column of $R$, taking every second column and ``raising" all columns appropriately to obtain a lower-triangular matrix, we obtain a matrix $B$. This matrix $B$ is then a Riordan array, given by
$$B=\left(g \frac{f}{x}, \frac{f^2}{x}\right).$$
\noindent
We have $$B = \left(\frac{f}{x}, x\right)\cdot A.$$
The general term $B_{n,k}$ of B is given by
\begin{eqnarray*}
B_{n,k}&=& [x^n] g\frac{f}{x}\left(\frac{f^2}{x}\right)^k \\
&=& [x^n] g x^{-k-1} f^{2k+1} \\
&=& [x^{n+k+1}] g f^{2k+1} \\
&=& T_{n+k+1,2k+1}.\end{eqnarray*}

\begin{example} We take the example of the binomial matrix
$$R=\left(\frac{1}{1-x}, \frac{x}{1-x}\right).$$
We then have
$$A=\left(\frac{1}{1-x}, \frac{x}{(1-x)^2}\right) \quad\textrm{with general term}\quad \binom{n+k}{2k},$$
and
$$B=\left(\frac{1}{(1-x)^2}, \frac{x}{(1-x)^2}\right) \quad\textrm{with general term}\quad \binom{n+k+1}{2k+1}.$$
\noindent The following decomposition makes this clear.
\begin{displaymath}\left(\begin{array}{ccccccc} \color{blue}{1} & 0 & 0 & 0
&0 & 0 & \cdots \\\color{blue}{1} & \color{red}{1} & 0 & 0 & 0 & 0 & \cdots \\ \color{blue}{1} & \color{red}{2} & \color{blue}{1} &
0 & 0 & 0 &
\cdots \\ \color{blue}{1} & \color{red}{3} & \color{blue}{3} & \color{red}{1} & 0 & 0 & \cdots \\ \color{blue}{1} & \color{red}{4} & \color{blue}{6} & \color{red}{4} &
\color{blue}{1} & 0 & \cdots \\\color{blue}{1} & \color{red}{5}  & \color{blue}{10} & \color{red}{10} & \color{blue}{5} & \color{red}{1} &\cdots\\ \vdots
& \vdots &
\vdots & \vdots & \vdots & \vdots &
\ddots\end{array}\right). \end{displaymath}
\noindent The matrices $A$ and $B$ are the coefficient arrays of the Morgan Voyce polynomials $b_n(x)$ and $B_n(x)$, respectively.
\end{example}

\begin{example} We take the Riordan array 
$$R=(c(x), xc(x))$$ where 
$$c(x)=\frac{1-\sqrt{1-4x}}{2x}.$$ 
Then we find that 
$$A=(c(x), xc(x)^2), \quad\quad B=(c(x)^2, xc(x)^2).$$ The matrix $R$ begins 
\begin{displaymath}\left(\begin{array}{ccccccc} \color{blue}{1} & 0 & 0 & 0
&0 & 0 & \cdots \\\color{blue}{1} & \color{red}{1} & 0 & 0 & 0 & 0 & \cdots \\ \color{blue}{2} & \color{red}{2} & \color{blue}{1} &
0 & 0 & 0 &
\cdots \\ \color{blue}{5} & \color{red}{5} & \color{blue}{3} & \color{red}{1} & 0 & 0 & \cdots \\ \color{blue}{14} & \color{red}{14} & \color{blue}{9} & \color{red}{4} &
\color{blue}{1} & 0 & \cdots \\\color{blue}{42} & \color{red}{42}  & \color{blue}{28} & \color{red}{14} & \color{blue}{5} & \color{red}{1} &\cdots\\ \vdots
& \vdots &
\vdots & \vdots & \vdots & \vdots &
\ddots\end{array}\right). \end{displaymath} 
We note that the matrix $A=(c(x), xc(x)^2)$ is the moment array for the family of orthogonal polynomials with coefficient array given by 
$$A^{-1}=(c(x), xc(x)^2)^{-1}=\left(\frac{1}{1+x}, \frac{x}{(1+x)^2}\right).$$ Denoting this family by $P_n(x)$, we have 
$$P_n(x)=(x-2)P_{n-1}(x)-P_{n-2}(x),$$ with $P_0(x)=1$, and $P_1(x)=x-1$. 
Similarly the matrix $B=(c(x)^2, xc(x)^2)$ is the moment array for the family of orthogonal polynomials with coefficient array given by $$B^{-1}=\left(\frac{1}{(1+x)^2}, \frac{x}{(1+x)^2}\right).$$ Denoting this family by $Q_n(x)$, we have 
$$Q_n(x)=(x-2)Q_{n-1}(x)-Q_{n-2}(x),$$ with $Q_0(x)=1$, and $Q_1(x)=x-2$. 

The inverse matrix $R^{-1}$ is given by 
\begin{displaymath}\left(\begin{array}{ccccccc} \color{red}{1} & 0 & 0 & 0
&0 & 0 & \cdots \\\color{blue}{-1} & \color{blue}{1} & 0 & 0 & 0 & 0 & \cdots \\ 0 & \color{red}{-2} & \color{red}{1} &
0 & 0 & 0 &
\cdots \\ 0 & \color{blue}{1} & \color{blue}{-3} & \color{blue}{1} & 0 & 0 & \cdots \\ 0 & 0 & \color{red}{3} & \color{red}{-4} &
\color{red}{1} & 0 & \cdots \\ 0 & 0  & \color{blue}{-1} & \color{blue}{6} & \color{blue}{-5} & \color{blue}{1} &\cdots\\ \vdots
& \vdots &
\vdots & \vdots & \vdots & \vdots &
\ddots\end{array}\right),  \end{displaymath} which is the Riordan array $\left(1-x, x(1-x)\right)$. In it we see the elements of $\color{blue}{A^{-1}}$ and $\color{red}{B^{-1}}$ in staggered fashion. 
\end{example}

\section{A counter-example}
It is natural to ask the question: is a matrix that contains two embedded Riordan arrays as above itself a Riordan array? The following example shows that this is not a sufficient condition on an array to be Riordan.
\begin{example} We shall construct an invertible integer lower-triangular matrix which has two embedded Riordan arrays in the fashion above,  but which is not itself a Riordan array.
We start with the essentially two-period sequence $(a_n)_{n \ge 0}$
$$1,2,3,2,3,2,3,\ldots.$$ We form the matrix
\begin{displaymath}\left(\begin{array}{ccccccc} 1 & 0 & 0 & 0
&0 & 0 & \cdots \\-2 & 1 & 0 & 0 & 0 & 0 & \cdots \\ 0 & -3 & 1 &
0 & 0 & 0 &
\cdots \\ 0 & 0 & -2 & 1 & 0 & 0 & \cdots \\ 0 & 0 & 0 & -3 &
1 & 0 & \cdots \\0 & 0  & 0 & 0 & -2 & 1 &\cdots\\ \vdots
& \vdots &
\vdots & \vdots & \vdots & \vdots &
\ddots\end{array}\right). \end{displaymath} The inverse of this matrix begins
\begin{displaymath}\left(\begin{array}{ccccccc} \color{blue}{1} & 0 & 0 & 0
&0 & 0 & \cdots \\\color{blue}{2} & \color{red}{1} & 0 & 0 & 0 & 0 & \cdots \\ \color{blue}{6} & \color{red}{3} & \color{blue}{1} &
0 & 0 & 0 &
\cdots \\ \color{blue}{12} & \color{red}{6} & \color{blue}{2} & \color{red}{1} & 0 & 0 & \cdots \\ \color{blue}{36} & \color{red}{18} & \color{blue}{6} & \color{red}{3} &
\color{blue}{1} & 0 & \cdots \\ \color{blue}{72} & \color{red}{36}  & \color{blue}{12} & \color{red}{6} & \color{blue}{2} & \color{red}{1} &\cdots\\ \vdots
& \vdots &
\vdots & \vdots & \vdots & \vdots &
\ddots\end{array}\right),\end{displaymath} where we note an alternating pattern of constant columns (with generating functions $\frac{1+2x}{1-6x^2}$ and $\frac{1+3x}{1-6x^2}$ respectively). Removing the first row of this matrix provides us with a production matrix, which is not of the form that produces a Riordan array (after the first column, subsequent columns would be shifted versions of the second column \cite{Prod1, Prod2}). Thus the resulting matrix will not be a Riordan array. This resulting matrix begins
\begin{displaymath}\left(\begin{array}{ccccccc} \color{blue}{1} & 0 & 0 & 0
&0 & 0 & \cdots \\\color{blue}{2} & \color{red}{1} & 0 & 0 & 0 & 0 & \cdots \\ \color{blue}{10} & \color{red}{5} & \color{blue}{1} &
0 & 0 & 0 &
\cdots \\ \color{blue}{62} & \color{red}{31} & \color{blue}{7} & \color{red}{1} & 0 & 0 & \cdots \\ \color{blue}{430} & \color{red}{215} & \color{blue}{51} & \color{red}{10} &
\color{blue}{1} & 0 & \cdots \\ \color{blue}{3194} & \color{red}{1597}  & \color{blue}{389} & \color{red}{87} & \color{blue}{12} & \color{red}{1} &\cdots\\ \vdots
& \vdots &
\vdots & \vdots & \vdots & \vdots &
\ddots\end{array}\right). \end{displaymath} \noindent
We now observe that for this matrix, we have
$$A=\left(\frac{1}{1+2x}, \frac{x}{1+5x+6x^2}\right)^{-1},$$
and
$$B=\left(\frac{1}{1+5x+6x^2}, \frac{x}{1+5x+6x^2}\right)^{-1}.$$
\noindent We notice that the sequence
$$1,2,10,62,430,3194,\ldots$$ has generating function given by the continued fraction
$$\cfrac{1}{1-\cfrac{2x}{1-\cfrac{3x}{1-\cfrac{2x}{1-\cdots}}}},$$ and secondly that
$$1+5x+6x^2=1+(2+3)x+2.3x^2=(1+2x)(1+3x).$$
\noindent This construction is easily generalized.

\end{example}
\section{Embedding a Riordan array}
Another natural question to ask is: if we are given a Riordan array $A$, is it possible to embed it as above into a Riordan array $R$? For this, we let
$$A=(u, v),$$ and seek to determine
$$R=(g, f)$$ such that $A$ embeds into $R$. For this, we need
$$ u=g, \quad \textrm{and} \quad v=\frac{f^2}{x}.$$ \noindent Thus we require that
$$f=\sqrt{xv}=x \sqrt{\frac{v}{x}}.$$
\noindent Since we are working in the context of integer valued Riordan arrays, we require that $v$ be such that $\sqrt{\frac{v}{x}}$ generates an integer sequence. We can state our result as follows.
\begin{proposition} The Riordan array
$$A=(u,v)$$ can be embedded in the Riordan array
$$R=\left(g, x \sqrt{\frac{v}{x}}\right)$$ on condition that $\sqrt{\frac{v}{x}}$ is the generating function of an integer sequence.
\end{proposition}
\begin{example} The Riordan array
$$A=\left(\frac{1}{\sqrt{1-4x}}, \frac{x}{1-4x}\right)$$ can be embedded in the Riordan array
$$R=\left(\frac{1}{\sqrt{1-4x}}, \frac{x}{\sqrt{1-4x}}\right).$$ \noindent
For this example, the matrix $A$ begins
\begin{displaymath}\left(\begin{array}{ccccccc} 1 & 0 & 0 & 0
&0 & 0 & \cdots \\2 & 1 & 0 & 0 & 0 & 0 & \cdots \\ 6 & 6 & 1 &
0 & 0 & 0 &
\cdots \\ 20 & 30 & 10 & 1 & 0 & 0 & \cdots \\ 70 & 140 & 70 & 14 &
1 & 0 & \cdots \\252 & 630  & 420 & 126 & 18 & 1 &\cdots\\ \vdots
& \vdots &
\vdots & \vdots & \vdots & \vdots &
\ddots\end{array}\right), \end{displaymath} while $R$ begins
\begin{displaymath}\left(\begin{array}{ccccccc} \color{blue}{1} & 0 & 0 & 0
&0 & 0 & \cdots \\\color{blue}{2} & \color{red}{1} & 0 & 0 & 0 & 0 & \cdots \\ \color{blue}{6} & \color{red}{4} & \color{blue}{1} &
0 & 0 & 0 &
\cdots \\ \color{blue}{20} & \color{red}{16} & \color{blue}{6} & \color{red}{1} & 0 & 0 & \cdots \\ \color{blue}{70} & \color{red}{64} & \color{blue}{30} & \color{red}{8} &
\color{blue}{1} & 0 & \cdots \\ \color{blue}{252} & \color{red}{256}  & \color{blue}{140} & \color{red}{48} & \color{blue}{10} & \color{red}{1} &\cdots\\ \vdots
& \vdots &
\vdots & \vdots & \vdots & \vdots &
\ddots\end{array}\right). \end{displaymath}
\end{example}

\section{A cascading decomposition}
We note that we can ``cascade" this embedding process, in the sense that given a Riordan array $R$, with embedded Riordan arrays $A$ and $B$, we can consider decomposing $A$ and $B$ in their turns and then continue this process. For instance, we can decompose
$$A=\left(g, \frac{f^2}{x}\right)$$ into the two matrices
$$A_A=\left(g, \frac{f^4}{x^3}\right), \quad \textrm{and} \quad B_A=\left(g\frac{f^2}{x^2}, \frac{f^4}{x^3}\right).$$ \noindent In their turn $A_A$ and $B_A$ can be decomposed and so on.

\section{Embeddings and orthogonal polynomials}
The phenomenon of embeddings as described above is not confined to Riordan arrays, as the continued fraction example above shows. To further amplify this point, we give another example involving a continued fraction. Although we take a particular case, the general case can be inferred easily from it. Thus we take the particular case of the continued fraction
$$\cfrac{1}{1-\cfrac{2x}{1-\cfrac{3x}{1-\cfrac{5x}{1-\cfrac{2x}{1-\cfrac{3x}{1-\cfrac{5x}{1-\cdots}}}}}}}.$$ This continued fraction is equal to
$$\cfrac{1}{1-2x-\cfrac{6x^2}{1-8x-\cfrac{10x^2}{1-5x-\cfrac{15x^2}{1-7x-\cfrac{6x^2}{1-8x-\cfrac{10x^2}{1-5x-\cdots}}}}}}.$$
By the theory of orthogonal polynomials, the power series expressed by both continued fractions is the generating function for the moment sequence of the family of orthogonal polynomials whose moment matrix (the inverse of the coefficient array of the orthogonal polynomials) has production matrix given by
\begin{displaymath}\left(\begin{array}{ccccccc} 2 & 1 & 0 & 0
&0 & 0 & \cdots \\6& 8 &1 & 0 & 0 & 0 & \cdots \\ 0 & 10 & 5 &
1 & 0 & 0 &
\cdots \\ 0 & 0 & 15 &7 & 1 & 0 & \cdots \\ 0 & 0 & 0 & 6 &
8 & 1 & \cdots \\0 & 0  & 0 & 0 & 10 & 5 &\cdots\\ \vdots
& \vdots &
\vdots & \vdots & \vdots & \vdots &
\ddots\end{array}\right). \end{displaymath} This production matrix generates the moment matrix $A$ that begins
\begin{displaymath}\left(\begin{array}{ccccccc} 1 & 0 & 0 & 0
&0 & 0 & \cdots \\2 & 1 & 0 & 0 & 0 & 0 & \cdots \\ 10 & 10 & 1 &
0 & 0 & 0 &
\cdots \\ 80 & 100 & 15 & 1 & 0 & 0 & \cdots \\ 760 & 1030 & 190 & 22 &
1 & 0 & \cdots \\7700 & 10900  & 2310 & 350 & 30 & 1 &\cdots\\ \vdots
& \vdots &
\vdots & \vdots & \vdots & \vdots &
\ddots\end{array}\right). \end{displaymath} In order to produce an embedding for this matrix, we proceed as follows. We form the matrix
\begin{displaymath}\left(\begin{array}{ccccccc} 1 & 0 & 0 & 0
&0 & 0 & \cdots \\-2 & 1 & 0 & 0 & 0 & 0 & \cdots \\ 0 & -3 & 1 &
0 & 0 & 0 &
\cdots \\ 0 & 0 & -5 & 1 & 0 & 0 & \cdots \\ 0 & 0 & 0 & -2 &
1 & 0 & \cdots \\0 & 0  & 0 & 0 & -3 & 1 &\cdots\\ \vdots
& \vdots &
\vdots & \vdots & \vdots & \vdots &
\ddots\end{array}\right). \end{displaymath} We invert this matrix, remove the first row of the resulting matrix, and use this new matrix as a production matrix. The generated matrix $R$ then begins
\begin{displaymath}\left(\begin{array}{ccccccc} 1 & 0 & 0 & 0
&0 & 0 & \cdots \\2 & 1 & 0 & 0 & 0 & 0 & \cdots \\ 10 & 5 & 1 &
0 & 0 & 0 &
\cdots \\ 80 & 40 & 10 & 1 & 0 & 0 & \cdots \\ 760 & 380 & 100 & 12 &
1 & 0 & \cdots \\7700 & 3850  & 1030 & 130 & 15 & 1 &\cdots\\ \vdots
& \vdots &
\vdots & \vdots & \vdots & \vdots &
\ddots\end{array}\right). \end{displaymath}
The moment matrix $A$ is evidently embedded in the matrix $R$. We can show that the corresponding matrix $B$ is the moment matrix for the family of orthogonal polynomials whose moments have generating function given by
$$\cfrac{1}{1-5x-\cfrac{15x^2}{1-7x-\cfrac{6x^2}{1-8x-\cfrac{10x^2}{1-5x-\cdots}}}}.$$ \noindent The matrix $B$ begins
\begin{displaymath}\left(\begin{array}{ccccccc} 1 & 0 & 0 & 0
&0 & 0 & \cdots \\5 & 1 & 0 & 0 & 0 & 0 & \cdots \\ 40 & 12 & 1 &
0 & 0 & 0 &
\cdots \\ 380 & 130 & 20 & 1 & 0 & 0 & \cdots \\ 3850 & 1410 & 300 & 25 &
1 & 0 & \cdots \\40400 & 15520  & 4060 & 440 & 32 & 1 &\cdots\\ \vdots
& \vdots &
\vdots & \vdots & \vdots & \vdots &
\ddots\end{array}\right), \end{displaymath} and it has production matrix
\begin{displaymath}\left(\begin{array}{ccccccc}
5 & 1 & 0 & 0 & 0 & 0 & \cdots \\
15& 7 & 1 & 0 & 0 & 0 & \cdots \\
0 & 6 & 8 & 1 & 0 & 0 & \cdots \\
0 & 0 & 10 & 5 & 1 & 0 & \cdots \\
0 & 0 & 0 & 15 & 7 & 1 & \cdots \\
0 & 0  & 0 & 0 & 6 & 8 &\cdots\\ \vdots
& \vdots &
\vdots & \vdots & \vdots & \vdots &
\ddots\end{array}\right). \end{displaymath}

The matrix $R^{-1}$ can now be characterized as the coefficient array of a family of polynomials $R_n(x)$ defined as follows. We let $P_n(x)$ be the family of orthogonal polynomials with coefficient array $A^{-1}$, and we let $Q_n(x)$ be the family of orthogonal polynomials with coefficient $B^{-1}$. Then we
have
\begin{displaymath}
R_n(x)=\begin{cases}
Q_{\frac{n}{2}}(x)x^{\frac{n}{2}},\quad \textrm{if $n$ is even;}\\
P_{\lceil \frac{n}{2}\rceil }(x)x^{\lfloor \frac{n}{2} \rfloor}, \quad \textrm{otherwise}.
\end{cases}\end{displaymath}

In the general case of a moment sequence generated by the continued fraction 
$$\cfrac{1}{1-\cfrac{\alpha x}{1-\cfrac{\beta x}{1-\cfrac{\gamma x}{1-\cfrac{\alpha x}{1-\cdots}}}}},$$ the matrix $A$ will be generated by the production matrix 
\begin{displaymath}\left(\begin{array}{ccccccc} \alpha & 1 & 0 & 0
&0 & 0 & \cdots \\ \alpha \beta & \beta+\gamma &1 & 0 & 0 & 0 & \cdots \\ 0 & \alpha \gamma & \alpha+\beta  &
1 & 0 & 0 &
\cdots \\ 0 & 0 & \beta \gamma  & \alpha+\gamma  & 1 & 0 & \cdots \\ 0 & 0 & 0 & \alpha \beta  &
\beta+\gamma  & 1 & \cdots \\0 & 0  & 0 & 0 & \alpha \gamma  & \alpha+\beta  &\cdots\\ \vdots
& \vdots &
\vdots & \vdots & \vdots & \vdots &
\ddots\end{array}\right), \end{displaymath} while the matrix $B$ is generated by the production matrix 
\begin{displaymath}\left(\begin{array}{ccccccc} 
\alpha+\beta & 1 & 0 & 0  &0 & 0 & \cdots \\
\beta \gamma & \alpha+\gamma &1 & 0 & 0 & 0 & \cdots \\ 
0 & \alpha \beta & \beta+\gamma &1 & 0 & 0 & \cdots \\ 
0 & 0 & \alpha \gamma &\alpha+\beta & 1 & 0 & \cdots \\ 
0 & 0 & 0 & \beta \gamma &  \alpha+\gamma & 1 & \cdots \\
0 & 0  & 0 & 0 & \alpha \beta & \beta +\gamma &\cdots\\ \vdots
& \vdots &
\vdots & \vdots & \vdots & \vdots &
\ddots\end{array}\right). \end{displaymath}
\section{Embeddings, orthogonal polynomials and Riordan arrays}
In this section, we consider the case of two related families of orthogonal polynomials $P_n(x)$ and $Q_n(x)$ defined by
$$P_n(x)=(x-7)P_{n-1}(x)-12P_{n-2}, $$ with
$P_0(x)=1$, $P_1(x)=x-3$, and
$$Q_n(x)=(x-7)Q_{n-1}(x)-12Q_{n-2}, $$ with
$Q_0(x)=1$, $Q_1(x)=x-7$.
We note that
$$(1+3x)(1+4x)=1+7x+12x^2.$$
The coefficient array of the polynomials $P_n(x)$ is then given by the Riordan array
$$A=\left(\frac{1}{1+3x}, \frac{x}{1+7x+12x^2}\right)^{-1},$$ while that of $Q_n(x)$ is given by
$$B=\left(\frac{1}{1+7x+12x^2}, \frac{x}{1+7x+12x^2}\right)^{-1}.$$
\noindent The matrix $A^{-1}$ begins
\begin{displaymath}\left(\begin{array}{ccccccc} 1 & 0 & 0 & 0
&0 & 0 & \cdots \\ 3 & 1 & 0 & 0 & 0 & 0 & \cdots \\ 21 & 10 & 1 &
0 & 0 & 0 &
\cdots \\ 183 & 103 & 17 & 1 & 0 & 0 & \cdots \\ 1785 & 1108 & 234 & 24 &
1 & 0 & \cdots \\18651 & 12349  & 3034 & 414 & 31 & 1 &\cdots\\ \vdots
& \vdots &
\vdots & \vdots & \vdots & \vdots &
\ddots\end{array}\right), \end{displaymath} while  $B^{-1}$ starts
\begin{displaymath}\left(\begin{array}{ccccccc} 1 & 0 & 0 & 0
&0 & 0 & \cdots \\ 7 & 1 & 0 & 0 & 0 & 0 & \cdots \\ 61 & 14 & 1 &
0 & 0 & 0 &
\cdots \\ 595 & 171 & 21 & 1 & 0 & 0 & \cdots \\ 6217 & 2044 & 330 & 28 &
1 & 0 & \cdots \\68047 & 24485  & 4690 & 538 & 35 & 1 &\cdots\\ \vdots
& \vdots &
\vdots & \vdots & \vdots & \vdots &
\ddots\end{array}\right). \end{displaymath} \noindent Letting
\begin{displaymath}
R_n(x)=\begin{cases}
Q_{\frac{n}{2}}(x)x^{\frac{n}{2}},\quad \textrm{if $n$ is even;}\\
P_{\frac{n+1}{2}}(x)x^{\lfloor \frac{n}{2} \rfloor}, \quad \textrm{otherwise},
\end{cases}\end{displaymath} we find that the inverse $R$ of the coefficient array of the family $R_n(x)$ is given by
\begin{displaymath}\left(\begin{array}{ccccccc} \color{blue}{1} & 0 & 0 & 0
&0 & 0 & \cdots \\\color{blue}{3} & \color{red}{1} & 0 & 0 & 0 & 0 & \cdots \\ \color{blue}{21} & \color{red}{7} & \color{blue}{1} &
0 & 0 & 0 &
\cdots \\ \color{blue}{183} & \color{red}{61} & \color{blue}{10} & \color{red}{1} & 0 & 0 & \cdots \\ \color{blue}{1785} & \color{red}{595} & \color{blue}{103} & \color{red}{14} &
\color{blue}{1} & 0 & \cdots \\ \color{blue}{18651} & \color{red}{6217}  & \color{blue}{1108} & \color{red}{171} & \color{blue}{17} & \color{red}{1} &\cdots\\ \vdots
& \vdots &
\vdots & \vdots & \vdots & \vdots &
\ddots\end{array}\right). \end{displaymath}
\noindent Thus the two Riordan arrays $A$ and $B$, which are the moment arrays of the two families of orthogonal polynomials $P_n(x)$ and $Q_n(x)$, respectively, embed into the generalized moment array $R$ for the family of polynomials $R_n(x)$.  Now the production matrix of $R$ begins
\begin{displaymath}\left(\begin{array}{ccccccc}
3 & 1 & 0 & 0 & 0 & 0 & \cdots \\
12& 4 & 1 & 0 & 0 & 0 & \cdots \\
36 & 12 & 3 & 1 & 0 & 0 & \cdots \\
144 & 48 & 12 & 4 & 1 & 0 & \cdots \\
432 & 144 & 36 & 12 & 3 & 1 & \cdots \\
1728 & 576  & 144 & 48 & 12 & 4 &\cdots\\ \vdots
& \vdots &
\vdots & \vdots & \vdots & \vdots &
\ddots\end{array}\right), \end{displaymath} where the columns have generating functions $\frac{1+4x}{1-12x^2}$, $\frac{1+3}{1-12x^2}$, respectively.
We now observe that this matrix is obtained by removing the first row of the inverse of the matrix
\begin{displaymath}\left(\begin{array}{ccccccc} 1 & 0 & 0 & 0
&0 & 0 & \cdots \\ -3 & 1 & 0 & 0 & 0 & 0 & \cdots \\ 0 & -4 & 1 &
0 & 0 & 0 &
\cdots \\ 0 & 0 & -3 & 1 & 0 & 0 & \cdots \\ 0 & 0 & 0 & -4 &
1 & 0 & \cdots \\

0 & 0  & 0 & 0 & -3 & 1 &\cdots\\ \vdots
& \vdots &
\vdots & \vdots & \vdots & \vdots &
\ddots\end{array}\right). \end{displaymath}
We note finally that the sequence
$$1, 3, 21, 183, 1785, 18651, 204141, \ldots$$ has generating function given by
$$g(x)=\cfrac{1}{1-\cfrac{3x}{1-\cfrac{4x}{1-\cfrac{3x}{1-\cdots}}}}.$$ We have in fact that
$$g(x)=\frac{1}{x}\textrm{Rev}\frac{x(1-4x)}{1-x}.$$
On the other hand, if we let 
$$P_n(x)=(x-7)P_{n-1}(x)-12 P_{n-2}(x),$$ but this time take $P_0(x)=1$ and 
$P_1(x)=x-4$, then the matrix $A$ becomes 
$$A=\left(\frac{1}{1+4x}, \frac{x}{1+7x+12x^2}\right)^{-1}.$$ The matrix $A$ then begins 
\begin{displaymath}\left(\begin{array}{ccccccc} 1 & 0 & 0 & 0
&0 & 0 & \cdots \\ 4 & 1 & 0 & 0 & 0 & 0 & \cdots \\ 28 & 11 & 1 &
0 & 0 & 0 &
\cdots \\ 244 & 117 & 18 & 1 & 0 & 0 & \cdots \\ 2380 & 1279 & 255 & 25 &
1 & 0 & \cdots \\

24868 & 14393  & 3364 & 442 & 32 & 1 &\cdots\\ \vdots
& \vdots &
\vdots & \vdots & \vdots & \vdots &
\ddots\end{array}\right), \end{displaymath} where the moment sequence 
$$1,4,28,244,2380,\ldots$$ has generating function 
$$\cfrac{1}{1-\cfrac{4x}{1-\cfrac{3x}{1-\cfrac{4x}{1-\cdots}}}},$$ or equivalently
$$\cfrac{1}{1-4x-\cfrac{12x^2}{1-7x-\cfrac{12x^2}{1-7x-\cdots}}}.$$ In this case, the matrix $R$ begins 
\begin{displaymath}\left(\begin{array}{ccccccc} \color{blue}{1} & 0 & 0 & 0
&0 & 0 & \cdots \\\color{blue}{4} & \color{red}{1} & 0 & 0 & 0 & 0 & \cdots \\ \color{blue}{28} & \color{red}{7} & \color{blue}{1} &
0 & 0 & 0 &
\cdots \\ \color{blue}{244} & \color{red}{61} & \color{blue}{11} & \color{red}{1} & 0 & 0 & \cdots \\ \color{blue}{2380} & \color{red}{595} & \color{blue}{117} & \color{red}{14} &
\color{blue}{1} & 0 & \cdots \\ \color{blue}{24868} & \color{red}{6217}  & \color{blue}{1279} & \color{red}{171} & \color{blue}{18} & \color{red}{1} &\cdots\\ \vdots
& \vdots &
\vdots & \vdots & \vdots & \vdots &
\ddots\end{array}\right). \end{displaymath} \noindent This matrix is then associated with the matrix

\begin{displaymath}\left(\begin{array}{ccccccc} 1 & 0 & 0 & 0
&0 & 0 & \cdots \\ -4 & 1 & 0 & 0 & 0 & 0 & \cdots \\ 0 & -3 & 1 &
0 & 0 & 0 &
\cdots \\ 0 & 0 & -4 & 1 & 0 & 0 & \cdots \\ 0 & 0 & 0 & -3 &
1 & 0 & \cdots \\

0 & 0  & 0 & 0 & -4 & 1 &\cdots\\ \vdots
& \vdots &
\vdots & \vdots & \vdots & \vdots &
\ddots\end{array}\right). \end{displaymath}

\bigskip
\hrule
\bigskip
\noindent 2010 {\it Mathematics Subject Classification}: Primary
15B36; Secondary 42C05, 33C45, 11B83.

\noindent \emph{Keywords:} Riordan array, production matrix, orthogonal polynomials, moment, integer sequence.


\begin{thebibliography}{13}


\bibitem{Barry_Meixner} P. Barry and A. Hennessy,\\Meixner-type results for Riordan arrays and associated
    integer sequences, \emph{J. Integer Seq.}, \textbf{13} (2010),
    \href{http://www.cs.uwaterloo.ca/journals/JIS/VOL13/Barry5/barry96s.pdf}{ Article 10.9.4}

\bibitem{Barry_Moment} P. Barry,
Riordan arrays, orthogonal polynomials as moments, and Hankel transforms,
\emph{J. Integer Seq.}, \textbf{14} (2011),
\href{http://www.cs.uwaterloo.ca/journals/JIS/VOL14/Barry1/barry97r2.html}{Article 11.2.2}

\bibitem{Cheon} G-S. Cheon, H. Kim, and L. W. Shapiro, Riordan group involution, \emph{Linear Algebra Appl.},
    \textbf{428}
    (2008), 941--952.

\bibitem{Chihara} T.~S. Chihara,  {\it An Introduction to Orthogonal Polynomials},
Gordon and Breach, New York, 1978.

\bibitem{CMS} C. Corsani, D. Merlini, and R. Sprugnoli,
    Left-inversion of combinatorial sums, \emph{Discrete
    Math.}
    \textbf{180} (1998), 107--122.

\bibitem{Prod1} E. Deutsch, L. Ferrari, and S. Rinaldi,
    Production matrices, \emph{Adv. in Appl.
    Math.},
    \textbf{34} (2005),
    101--122.

\bibitem{Prod2} E. Deutsch, L. Ferrari, and S. Rinaldi,
Production matrices and Riordan arrays, \emph{Ann. Comb.}, \textbf{13} (2009), 65--85.

\bibitem{Gautschi} W. Gautschi, {\it Orthogonal Polynomials: Computation and
Approximation}, Clarendon Press, Oxford.

\bibitem{He}
Tian-Xiao He, R. Sprugnoli, Sequence characterization of Riordan arrays, \emph{Discrete Math.} \textbf{2009} (2009),
 3962--3974.

\bibitem{Jin}
S.-T. Jin, A characterization of the Riordan Bell subgroup by C-sequences, \emph{Korean J. Math.} \textbf{17} (2009), 147--154.

\bibitem{Alter} D. Merlini, D.~G.~Rogers, R. Sprugnoli, and M.~C.~Verri,
    On some alternative characterizations of Riordan
    arrays,
    \emph{Canad. J. Math.}, \textbf{49} (1997),
     301--320.

\bibitem{SGWW} L.~W.~Shapiro, S. Getu, W.-J. Woan, and L.C. Woodson,
The Riordan Group, \emph{Discr. Appl. Math.} \textbf{34} (1991), 229--239.

\bibitem{SL1} N.~J.~A.~Sloane, \emph{The
On-Line Encyclopedia of Integer Sequences}. Published electronically
at \href{http:oeis.org}{\texttt{http://oeis.org}}, 2011.

\bibitem{SL2} N.~J.~A.~Sloane, The On-Line Encyclopedia of Integer
Sequences, \emph{Notices Amer. Math. Soc.}, \textbf{50} (2003),  912--915.

\bibitem{Spru} R. Sprugnoli, Riordan arrays and combinatorial sums,
\emph{Discrete Math.} \textbf{132} (1994), 267--290.

\bibitem{Szego} G. Szeg\"o, \emph{Orthogonal Polynomials}, 4th
    ed. Providence, RI, Amer. Math. Soc., 1975.

\bibitem{Wall} H.~S. Wall, \emph{Analytic Theory of
    Continued Fractions}, AMS Chelsea Publishing, 1967.


\end{thebibliography}
\end{document}